\documentclass[11pt,reqno]{amsart}
\usepackage{amsmath,amsfonts, amssymb, amsthm}
  \usepackage{paralist}
  \usepackage{graphics} 
  \usepackage{epsfig} 
\usepackage{graphicx}  \usepackage{epstopdf}
 \usepackage[colorlinks=true]{hyperref}
\hypersetup{urlcolor=blue, citecolor=red}

\newtheorem{theorem}{Theorem}[section]

\newtheorem{lemma}[theorem]{Lemma}

\newtheorem{conjecture}{Conjecture}[section]

\theoremstyle{definition}

\newtheorem{remark}{Remark}[section]

\def\pmod #1{\ ({\rm{mod}}\ #1)}
\def\Z{\Bbb Z}
\def\N{\Bbb N}

\def\Q{\Bbb Q}

\def\l{\left}
\def\r{\right}
\def\bg{\bigg}
\def\({\bg(}
\def\){\bg)}
\def\t{\text}
\def\f{\frac}
\def\mo{{\rm{mod}\ }}
\def\pmod#1{\ (\mo\ #1)}

\def\gs{\geq}
\def\se {\subseteq}
\def\sm{\setminus}

\def\al{\alpha}

\def\ve{\varepsilon}

\def\eq{\equiv}

\def\da{\delta}

\def\Proof{\noindent{\it Proof}}

\begin{document}
\hbox{Preprint}
\medskip

\title[On exponential diophantine equations over $\Q$ with few unknowns]
      {On exponential diophantine equations over $\Q$ with few unknowns}
\author[Zhi-Wei Sun]{Zhi-Wei Sun}

\address{Department of Mathematics, Nanjing
University, Nanjing 210093, People's Republic of China}
\email{zwsun@nju.edu.cn}

\keywords{Undecidability, exponential diophantine equations, recursively enumerable sets, Hilbert's Tenth Problem over $\Q$.
\newline \indent 2020 {\it Mathematics Subject Classification}. Primary 11D61, 03D35; Secondary
03D25, 11U05.
\newline \indent Supported by the Natural Science Foundation of China (grant no. 11971222).}

\begin{abstract} In this paper we obtain three undecidable results
for exponential diophantine equations over the field $\Q$ of rational numbers.
For example, we prove that there is no algorithm to decide the solvability of a general
exponential diophantine equation $F(x_1,\ldots,x_{8})=0$ over $\Q$ with eight unknowns.
\end{abstract}
\maketitle

\section{Introduction}

Exponential diophantine equations have the form
$$F_1(x_1,\ldots,x_n)-F_2(x_1,\ldots,x_n)=0,$$
where $F_1$ and $F_2$ are expressions constructed from variables and particular natural numbers using
addition, multiplication, and exponentiation. (Note that $0^0$ is regarded as $1$.)
Here is an example of exponential diophantine equation:
$$x^{2^{y^x}}+y^{x+3y}-(5z^{2x^2}+xyz+4)=0.$$
In 1961 M. Davis, H. Putnam and J. Robinson \cite{DPR}
proved that the solvability of of a general exponential diophantine equation over $\N$
is undecidable, this is the first major step towards Y. Matiyasevich's negative solution \cite{M70}
of Hilbert's Tenth Problem (HTP). In this direction, the 11 Unknowns Theorem established by Z.-W. Sun \cite{S21}
states that there is no algorithm to decide
 for any given polynomial $P(x_1,\ldots,x_{11})\in\Z[x_1,\ldots,x_{11}]$
whether the equation $P(x_1,\ldots,x_{11})=0$ has integer solutions.

HTP over $\Q$ (the field of rational numbers) is still open.
We may also consider the decision problem for exponential diophantine equations over $\Q$.
To avoid uncertainty like $(-1)^{2/3}$ ,
in exponential diophantine equations we use the exponential function $x^y$ only in the case $x,y\gs0$.
Recently, M. Prunescu \cite{Pr}
used a clever trick to show that there is no algorithm to decide whether an arbitrary exponential diophantine equation over $\Q$ has rational solutions.

In this paper, we establish the following three new theorems for exponential diophantine equations over $\Q$.

\begin{theorem}\label{Th-8} {\rm (i)} For any r.e. (recursively enumerable) subset $A$
of $\N=\{0,1,2,\ldots\}$, there is an exponential diophantine equation
$$F(x,x_1,\ldots,x_{8})=0$$
such that $a\in\N$ belongs to $A$ if and only if
$$F(a,x_1,\ldots,x_8)=0$$
for some $x_1,\ldots,x_8\in\Q$.

{\rm (ii)} There is no algorithm to decide for any exponential diophantine equation
$$F(x_1,\ldots,x_8)=0$$
with eight unknowns, whether it has solutions over $\Q$.
\end{theorem}

\begin{theorem}\label{Th-20} {\rm (i)} For any r.e. set $A\se \N$, there is an exponential diophantine equation
$$F(x,x_1,\ldots,x_{10})=0$$
such that $a\in\N$ belongs to $A$ if and only if
$$F(a,x_1^2,\ldots,x_{10}^2)=0$$
for some $x_1,\ldots,x_{10}\in\Q$.

{\rm (ii)} There is no algorithm to decide for any exponential diophantine equation
$$F(x_1,\ldots,x_{10})=0$$
with ten unknowns, whether the equation
\begin{equation}F(x_1^2,\ldots,x_{10}^2)=0
\end{equation} has solutions over $\Q$.
\end{theorem}
\begin{remark} In contrast with Theorem \ref{Th-20}(ii), Sun \cite{S21} proved that
there is no algorithm to decide
 for any given polynomial $P(x_1,\ldots,x_{17})\in\Z[x_1,\ldots,x_{17}]$
whether the equation $P(x_1^2,\ldots,x_{17}^2)=0$ has integer solutions.
\end{remark}

 The number $10$ in Theorem \ref{Th-20} seems not optimal, however it looks
challenging to replace it by a smaller number. Nevertheless, we pose the following conjecture.

\begin{conjecture} There is no algorithm to decide for any exponential diophantine equation
$F(x,y,z)=0$, whether $F(x^2,y^2,z^2)=0$ for some $x,y,z\in\Q$.
\end{conjecture}

\begin{theorem}\label{Th-11} Let $p_1,\ldots,p_{10}$ be ten distinct primes.

{\rm (i)} For any r.e. set $A\se \N$, there is a polynomial $P(u,v,x_1,\ldots,x_{11})$ with integer coefficients
such that $a\in\N$ belongs to $A$ if and only if
$$P\l(a,x_0,x_1,\ldots,x_{10},p_1^{x_1^2}\cdots p_{10}^{x_{10}^2}\r)=0$$
for some $x_0,\ldots,x_{10}\in\Q$.

{\rm (ii)} There is no algorithm to decide for any polynomial $P(x_0,\ldots,x_{11})$
with integer coefficients whether the equation
$$P\l(x_0,x_1,\ldots,x_{10},p_1^{x_1^2}\cdots p_{10}^{x_{10}^2}\r)=0$$
has solutions over $\Q$.
\end{theorem}

We will provide some lemmas in the next section, and prove Theorems \ref{Th-8}-\ref{Th-11} in Section 3.
All variables in Sections 2 and 3 range over $\Q$ unless specified.

Throughout this paper, we adopt logical symbols $\land$ (conjunction) and $\lor$ (disjunction), and also set
\begin{equation}\square:=\{r^2:\ r\in\Q\}.\end{equation}

\section{Some lemmas}
 \setcounter{equation}{0}
 \setcounter{conjecture}{0}
 \setcounter{theorem}{0}

 Prunescu \cite{Pr} showed that a rational number $w$ belongs to $\N$
 if and only if there are rational numbers $x>0$ and $y>x$ such that $x^y=y^x$ and $wy=(w+1)x$.
 This follows from the known result that the only rational solutions of the equation $x^y=y^x$
 with $0<x<y$ are given by
$$x=\l(1+\f1n\r)^n\ \ \t{and}\ \ y=\l(1+\f1n\r)^{n+1}\ \ (n\in\Z^+=\{1,2,3,\ldots\}),$$
 which dates back to L. Euler (cf. L. E. Dickson \cite[p.\,687]{Di}) and appeared in
 a modern reference \cite{Sv} with a detailed proof.

 Our following lemma gives a rather simple way to characterize $\Z$ in $\Q$.

\begin{lemma}\label{Au} Let $\al_1,\ldots,\al_k\in\Q$, and let $p_1,\ldots,p_k$ be distinct primes.
Then
\begin{equation} \label{S} \al_1,\ldots,\al_k\in\Z\iff \prod_{i=1}^kp_i^{\al_i}\in\Q
\iff \prod_{i=1}^kp_i^{\al_i^2}\in\Q.
\end{equation}
\end{lemma}
\Proof. For any $\al\in\Q$, if $\al^2$ is an integer $m$ then $\al$ is a rational algebraic integer
and hence $\al\in\Z$ by Prop. 6.1.1 of \cite[p.\,66]{IR}.
So we only need to show the first equivalence in \eqref{S}.

$\Rightarrow$: This direction is obvious.

$\Leftarrow$: Suppose that $\prod_{i=1}^kp_i^{\al_i}$ is a rational number $r$.
Choose $n\in\Z^+$ such that $m_i=n\al_i\in\Z$ for all $i=1,\ldots,k$. Then
\begin{equation}\label{r^n}\prod_{i=1}^kp_i^{m_i}=r^n.\end{equation}
For any prime $p$ and $x\in\Q\sm\{0\}$, we write $\nu_p(x)$ to denote the $p$-adic valuation of $x$.
By \eqref{r^n}, for any $i=1,\ldots,k$ we have
$$m_i=\nu_{p_i}(r^n)=n\nu_{p_i}(r)\eq0\pmod n.$$
Therefore $\al_i=m_i/n\in\Z$ for all $i=1,\ldots,k$.

In view of the above, we have proved Lemma \ref{Au}. \qed

\begin{remark} Let $a>1$ be an integer which is not a perfect power.
For a rational number $\al=m/n$ with $m\in\Z$, $n\in\Z^+$ and $\gcd(m,n)=1$,
if $a^{\al}=r\in\Q$ then $a^m=r^n$ and hence $a$ is an $n$-th power
(since $\gcd(m,n)=1$), so $a^\al\in\Q$ if and only if $\al\in\Z$.
We also note that a rational number $x\gs0$ is an integer if and only if $x^x\in\Q$. In fact, if $x=m/n$ with $m,n\in\Z^+$, $\gcd(m,n)=1$ and $n>1$, and $x^x=r\in\Q$, then
$(m/n)^m=r^n$, hence for any prime divisor $p$ of $n$ we have $n\mid \nu_p(n)$ and thus $n\gs p^n>n$
which is impossible.
\end{remark}

\begin{lemma}\label{4m+2} For any integer $m$, we have
\begin{equation} m\gs0\iff\exists x\in\Z[x\not=0\land (4m+2)x^2+1\in\square].
\end{equation}
\end{lemma}
\Proof. For any $x\in\Z$, the integer $(4m+2)x^2+1$ lies in $\square=\{r^2:\ r\in\Q\}$
if and only if it is an integer square. So, this lemma is essentially equivalent
to \cite[Lemma 2]{S92b} obtained via Pell equations. \qed

The following result is an analogue of Matiyasevich-Robinson's Relation-Combining Theorem
\cite{MR}.

\begin{lemma} [G.-R. Zhang and Z.-W. Sun \cite{ZS}]\label{Lem-ZS} Let $A_1,\ldots,A_k\in\Q\sm\{0\}$, and define
\begin{align*}{\mathcal J}_k(A_1,\ldots,A_k,x)=&\prod_{s=1}^kA_s^{(k-1)2^{k+1}}
\times\prod_{\ve_1,\ldots\ve_k\in\{\pm1\}}
\l(x+\sum_{s=1}^k\ve_s\sqrt{A_s}\,W^{s-1}\r),
\end{align*}
where
$$W=\l(k+\sum_{s=1}^kA_s^2\r)\l(1+\sum_{s=1}^kA_s^{-2}\r).$$
Then ${\mathcal J}_k(x_1,\ldots,x_k,x)$ is a polynomial with integer coefficients. Moreover,
\begin{equation}A_1,\ldots,A_k\in\square \iff \exists x [{\mathcal J}_k(A_1,\ldots,A_k,x)=0].
\end{equation}
\end{lemma}

\begin{lemma}[Matiyasevich, 1979]\label{Lem-M}
For any r.e. set $A\se\N$, there is an exponential diophantine equation
$$f(t,x,y,z)=0$$
such that for any $a\in\N$ we have
\begin{equation}\label{A}
a\in A\iff \exists x\in\N\,\exists y\in\N\,\exists z\in\N\,[f(a,x,y,z)=0].
\end{equation}
\end{lemma}
\begin{remark} This result of Matiyasevich \cite{M79} (see also Section 8.2 of \cite[pp.\,156--160]{M-book}) improves the Davis-Putnam-Robinson Theorem
\cite{DPR} greatly.
\end{remark}

\begin{lemma} \label{>=0} Let $\al$ be a rational number. Then
\begin{equation}\al\gs0\iff\exists x_1\exists x_2\exists x_3[\al=x_1^2+x_2^2+x_3^2\lor \al=x_1^2+x_2^2+2x_3^2].\end{equation}
\end{lemma}
\Proof. $\Leftarrow$: This is obvious.

$\Rightarrow$: Write $\al=a/b$ with $a\in\N$ and $b\in\Z^+$. By the theory of ternary quadratic forms,
$$\N\sm\{x^2+y^2+z^2:\ x,y,z\in\N\}=\{4^k(8m+7):\ k,m\in\N\}$$
and $$\N\sm\{x^2+y^2+2z^2:\ x,y,z\in\N\}=\{4^k(16m+14):\ k,m\in\N\}$$
(cf. \cite[pp.\,112--113]{D39}). Note that
$$\{4^k(8m+7):\ k,m\in\N\}\cap \{4^k(16m+14):\ k,m\in\N\}=\emptyset.$$
So, for some $\da\in\{1,2\}$ we have
$ab=x^2+y^2+\da z^2$ for some $x,y,z\in\N$. Hence
$$\al=\f{ab}{b^2}=\l(\f xb\r)^2+\l(\f yb\r)^2+\da\l(\f zb\r)^2.$$
This concludes the proof. \qed

\section{Proofs of Theorems \ref{Th-8}-\ref{Th-11}}
 \setcounter{equation}{0}
 \setcounter{conjecture}{0}
 \setcounter{theorem}{0}

 It is known that there are nonrecursive r.e. subsets of $\N$ (see, e.g., N. Cutland \cite[pp.\,140--141]{C}).
Thus, for each of Theorems \ref{Th-8}-\ref{Th-11}, its first part implies the second part. So it remains to prove the first parts of Theorems \ref{Th-8}-\ref{Th-11}.

 Let $A$ be any r.e. subset of $\N$. By Lemma \ref{Lem-M},
there is an exponential diophantine equation
$f(t,x,y,z)=0$ such that \eqref{A} holds for all $a\in \N$.

Let $a\in\N$. In view of Lemmas \ref{Au}-\ref{Lem-ZS}, we have
\begin{align*}&\exists x\in\N\,\exists y\in\N\,\exists z\in\N [f(a,x,y,z)=0]
\\\iff&\exists x\exists y\exists z\exists \bar x\exists\bar y\exists \bar z
[x,y,z,\bar x,\bar y,\bar z\in\Z\land \bar x\bar y\bar z\not=0
\land (4x+2)\bar x^2+1\in\square
\\& \land (4y+2)\bar y^2+1\in\square
\land  (4z+2)\bar z^2+1\in\square\land f(a,x,y,z)=0]
\\\iff&\exists x\exists y\exists z\exists \bar x\exists\bar y\exists \bar z\exists u\exists v
\big[\big(u\bar x\bar y\bar z 2^{x^2}3^{y^2}5^{z^2}7^{\bar x^2}11^{\bar y^2}13^{\bar z^2}-1\big)^2+f(a,x,y,z)^2
\\&+\mathcal{J}_3((4x+2)\bar x^2+1,(4y+2)\bar y^2+1,(4z+2)\bar z^2+1,v)^2=0\big].
\end{align*}
In light of Lemmas \ref{Au} and \ref{>=0}, we also have
\begin{align*}&\exists x\in\N\,\exists y\in\N\,\exists z\in\N [f(a,x,y,z)=0]
\\\iff& \exists x_1\exists x_2\exists x_3\exists y_1\exists y_2\exists y_3
\exists z_1\exists z_2\exists z_3
\exists\da_1\in\{1,2\}\exists\da_2\in\{1,2\}
\exists\da_3\in\{1,2\}
\\&[2^{x_1^2+x_2^2+\da_1 x_3^2}3^{y_1^2+y_2^2+\da_2 y_3^2}5^{z_1^2+z_2^2+{\da_3}z_3^2}\in\Q
\\&\land f(a,x_1^2+x_2^2+{\da_1}x_3^2,y_1^2+y_2^2+{\da_2}y_3^2,z_1^2+z_2^2+{\da_3}z_3^2)=0]
\\\iff&\exists wx_1\exists x_2\exists x_3\exists y_1\exists y_2\exists y_3
\exists z_1\exists z_2\exists z_3
\\&[F(a,w^2,x_1^2,x_2^2,x_3^2,y_1^2,y_2^2,y_3^2,z_1^2,z_2^2,z_3^2)=0],
\end{align*}
where $F(a,w^2,x_1^2,x_2^2,x_3^2,y_1^2,y_2^2,y_3^2,z_1^2,z_2^2,z_3^2)$ is the product of
those
\begin{align*}
&\l(w^2-(2^{x_1^2+x_2^2+{\da_1}x_3^2}3^{y_1^2+y_2^2+{\da_2}y_3^2}5^{z_1^2+z_2^2+{\da_3}z_3^2})^2\r)^2
\\&+f(a,x_1^2+x_2^2+{\da_1}x_3^2,y_1^2+y_2^2+{\da_2}y_3^2,
z_1^2+z_2^2+{\da_3}z_3^2)^2
\end{align*}
with $\da_1,\da_2,\da_3\in\{1,2\}$.

By Sun \cite[Theorem 1.1(ii)]{S21}, there is a polynomial $Q(x_0,x_1,\ldots,x_{10})\in\Z[x_0,\ldots,x_{10}]$
such that for any $a\in A$ we have
$$a\in A\iff\exists x_1\ldots\exists x_{10}[x_1,\ldots,x_{10}\in\Z\land x_{10}\not=0
\land Q(a,x_1,\ldots,x_{10})=0].$$
Combining this with Lemma \ref{Au}, for any $a\in\N$ we get
\begin{align*}a\in A&\iff \exists x_0\exists x_1\ldots\exists x_{10}
\big[\big(x_0x_{10}p_1^{x_1^2}\cdots p_{10}^{x_{10}^2}-1\big)^2+Q(a,x_1,\ldots,x_{10})^2=0\big].
\end{align*}

In view of the above, we have completed the proofs of Theorems \ref{Th-8}-\ref{Th-11}.



\begin{thebibliography}{99}

\bibitem {C} N. Cutland,  Computability, Cambridge Univ. Press, Cambridge, 1980.


\bibitem {DPR} M. Davis, H. Putnam  and J. Robinson,
 {\it  The decision problem for exponential diophantine equations},
  Ann. of Math. {\bf 74}(2) (1961), 425--436.

 \bibitem{D39} L. E. Dickson, Modern Elementary Theory of Numbers, University of Chicago Press, Chicago, 1939.

 \bibitem{Di} L. E. Dickson, History of the Theory of Numbers, Vol. II, AMS Chelsea Publ., 1999.

\bibitem{IR} K. Ireland and M. Rosen, A Classical Introduction to Modern Number Theory, 2nd Edition, Grad. Texts. Math., vol. 84, Springer, New York, 1990.

\bibitem {M70} Y. Matiyasevich, {\it Enumerable sets are diophantine},
 Dokl. Akad. Nauk SSSR {\bf 191} (1970), 279--282;
English translation with addendum, Soviet Math. Doklady {\bf 11} (1970), 354--357.

\bibitem{M79} Y. Matiyasevich, {\it Algorithmic unsolvability of exponential Diophantine equations
in three unknowns}, Slecta Math. Sovietica {\bf 3}(3) (1983/84), 223--232.

\bibitem{M-book} Y. Matiyasevich, Hilbert's Tenth Problem, MIT Press, Cambridge, Massachusetts,
1993.

\bibitem {MR} Y. Matiyasevich and J. Robinson,
 {\it Reduction of an arbitrary diophantine equation to one in 13 unknowns}, Acta Arith. {\bf 27} (1975), 521--553.

\bibitem{Pr} M. Prunescu, {\it The exponetial diophantine problem for $\Q$}, J. Symb. Log.
{\bf 85} (2020), 671--672.

\bibitem {S92b} Z.-W. Sun, {\it A new relation-combining theorem and its application},
 Z. Math. Logik Grundlag. Math. {\bf 38} (1992), 209--212.

\bibitem{S21} Z.-W. Sun, {\it Further results on Hilbert's Tenth Problem}, Sci. China Math.
{\bf 64} (2021), 281--306.

\bibitem{Sv} M. Sved, {\it On the rational solutions of $x^y=y^x$}, Math. Magzine {\bf 63} (1990), 30--33.

\bibitem {ZS} G.-R. Zhang and Z.-W. Sun, {\it $\Q\sm\Z$ is diophantine over $\Q$ with $32$ unknowns}, preprint, arXiv:2104.02520, 2021.

\end{thebibliography}
 \end{document}